\input ./style/arxiv-ba.cfg
\documentclass[ba,linksfromyear,preprint]{imsart}
\makeatletter
   \@ifpackageloaded{natbib}{}{\usepackage{natbib}}
\makeatother
\newcommand{\enquote}[1]{``#1''}
\defcitealias{FG2015}{FG}

\pubyear{2015}
\volume{10}
\issue{3}
\firstpage{741}
\lastpage{748}
\doi{10.1214/15-BA944B}% Updated by VTEXPTS2LaTeX.exe, 15.04.2015 16:02
\begin{document}

\begin{frontmatter}
\title{Comment~on~Article~by~Ferreira~and~Gamerman\thanksref{T1}}
%\thankstext{<label>}{<text>}
%
% Discussant paper:
%\relateddois{<label>}{Related articles:
%DOI:~\relateddoi[ms=BAXXX]{Related item:}{10.1214/00-BAXXX},
%DOI:~\relateddoi[ms=BAXXX]{Related item:}{10.1214/00-BAXXX},
%DOI:~\relateddoi[ms=BAXXX]{Related item:}{10.1214/00-BAXXX}; rejoinder
%at
%DOI:~\relateddoi[ms=BAXXX]{Related item:}{10.1214/00-BAXXX}.}
%
% Comment, Discussion, Rejoinder
%\relateddois{<label>}{Main article DOI: \relateddoi[ms=BAXXX]{Related
%item:}{10.1214/00-BAXXX}.}
%
\relateddois{T1}{Main article DOI: \relateddoi[ms=BA944]{Related item:}{10.1214/15-BA944}.}
\runtitle{Comment on Article by Ferreira and Gamerman}

\begin{aug}
\author[addr1]{\fnms{Noel} \snm{Cressie}\corref{}\ead[label=NCemail]{ncressie@uow.edu.au}}\thanksref{NCthanks}
\and
\author[addr2]{\fnms{Raymond L.} \snm{Chambers}}

\runauthor{N. Cressie and R. L. Chambers}

\address[addr1]{National Institute of Applied Statistics Research
Australia, University of Wollongong, Australia,\\ \printead{NCemail}}
\address[addr2]{National Institute of Applied Statistics Research
Australia, University of Wollongong, Australia}

\thankstext{NCthanks}{Cressie's research was partially supported by
the US National Science Foundation and the US Census Bureau through
the NSF-Census Research Network program; and it was partially supported
by a 2015-2017 Australian Research Council Discovery Project.}

\end{aug}

%% Abstract %%
%
\begin{abstract}
A utility-function approach to optimal spatial sampling design
is a powerful way to quantify what ``optimality'' means. The
emphasis then should be to capture all possible contributions
to utility, including scientific impact and the cost of sampling.
The resulting sampling plan should contain a component of designed
randomness that would allow for a non-parametric design-based analysis
if model-based assumptions were in doubt.
\end{abstract}

%% Keywords %%
%

% \begin{keyword}[class=MSC]
% \kwd[Primary ]{}
% \kwd[; secondary ]{}
% \end{keyword}

\end{frontmatter}

%% Mainmatter %%

\section{Introduction}\label{section:intro}

We would like to express our appreciation to Gustavo da Silva Ferreira
and Dani Gamerman (hereafter,~\citetalias{FG2015}) for their paper on
Bayesian preferential spatial sampling (\citealp{FG2015}) and to the
editor of \textit{Bayesian Analysis} for the opportunity to contribute
to the discussion. Building on the papers by~\cite{Muller1999} and~\cite
{Diggle2010}, the authors give a Bayesian approach to choosing \textit
{new} sampling locations after initial data are assumed to have been
obtained under preferential sampling.

\section{What is Fixed and What is Random?}
\label{sect:FixedRandom}

Let the initial sample be $\mathbf{y_x}$, obtained at preferential
sampling locations $\mathbf{x}$; note that we have emphasised
dependence of $\mathbf{y}$ on $\mathbf{x}$ through the notation $\mathbf
{y_x}$, but it is exactly the same as what~\citetalias{FG2015} notate
as $\mathbf{y}$. The observation locations $\mathbf{x}$ and the
observations $\mathbf{y_x}$ are known by the statistician designing the
next phase of the study, and hence all criteria and inferences should
depend on both $\mathbf{x}$ and $\mathbf{y_x}$. One can see this most
clearly in ~\citetalias{FG2015}'s definition of the Bayesian design
criterion $U(\mathbf{d})$, given at the beginning of~\citetalias
{FG2015}-Section 4. However, the reader should notice that
equations~\citetalias{FG2015}-(3) and~\citetalias{FG2015}-(4) do not
emphasise conditioning on $\mathbf{x}$, along with $\mathbf{y_x}$,
something we assume is an oversight on the part of the authors.

It helped us to augment the notation for the utility function from
$u(\mathbf{d},\theta,\mathbf{y_d})$ to $u(\mathbf{d},\theta,\mathbf
{y_d};\mathbf{x},\mathbf{y_x})$; and likewise we suggest that the
expected utility be notated:
\begin{equation}
\label{eq:exp-utility}
U(\mathbf{d};\mathbf{x},\mathbf{y_x})=E(u(\mathbf{d},\theta,\mathbf
{y_d};\mathbf{x},\mathbf{y_x})|\mathbf{x},\mathbf{y_x}),
\end{equation}
where $\mathbf{d}$ is considered fixed and the expectation is taken
over $[\theta,\mathbf{y_d}|\mathbf{x},\mathbf{y_x}]$. When there are
many ``stakeholders'' (e.g., in an environmental study), each coming
with his/her own utility function, how can a single utility function be
constructed?~\citet[Chapter 11]{Le2006} opt for one based on entropy.
Do the authors have any other suggestions to build ``compromise'' into
a utility function?

Notation is really important in these complex situations, so in the
case of the utility function defined by~\citetalias{FG2015}-(4), which
involves the latent process (not the observations), we suggest that $u$
be rewritten as:
\begin{equation*}
\label{eq:rewriting-u}
u(\mathbf{d},\theta,\mathbf{s_d};\mathbf{x},\mathbf{y_x});
\end{equation*}
that is, $\mathbf{s_d}$ replaces $\mathbf{y_d}$ in $u$. Depending on
the context, $u$ could be a function of the new observations, $\mathbf
{y_d}\equiv\left(y(d_1),\ldots,y(d_m)\right)^{\prime}$, or of the
corresponding latent process, $\mathbf{s_d}\equiv\left(s(d_1),\ldots
,s(d_m)\right)^{\prime}$; recall that~\citetalias{FG2015} have defined
$y(\cdot)$ to be a noisy, shifted version of the mean-zero latent
process $s(\cdot)$.

In the rest of our discussion, we follow the authors' lead and use
\eqref{eq:exp-utility}, albeit with our modified notation that
emphasises dependence on $\mathbf{x}$ and $\mathbf{y_x}$. The
utility-function approach to optimal design is attractive, but it will
only be truly useful when components that quantify ``how much?'' and
``why?'' are specifically included; see Sections~\ref{section:spatial}
and~\ref{section:utility} for further discussion.

As~\citetalias{FG2015} make clear, the process $s(\cdot)$, the sampling
locations $\mathbf{x}\equiv(x_1,\ldots,x_n)^{\prime}$, and the
observations $\mathbf{y_x}\equiv\left(y(x_1),\ldots,y(x_n)\right
)^{\prime}$ have a possibly complex joint distribution. Following~\cite
{Diggle2010}, the authors put structure on this joint distribution by
assuming~\citetalias{FG2015}-(1),~\citetalias{FG2015}-(2), and a
log-Gaussian Cox process. From the point of view of sample survey
design, the information in $\mathbf{x}$ and $\mathbf{y_x}$ is
comparable to what one would gain from a pilot study, but it requires
knowledge of components of $\theta$ in order to make the pilot study
operational. The following suggestion seems compatible with the
authors' approach to optimal design via preferential sampling.

It is hard to design a study if there is no knowledge from which to
draw. In the pre-pilot phase, one might choose a simple random sample
which, in the spatial context, means that observation locations are
sampled uniformly from the spatial region of interest, $A$. We note
that this corresponds to a degenerate case of preferential sampling
where $\beta$, the coefficient of $s(\cdot)$ in the log-intensity, is
equal to $0$, and we also note the presence of randomness in this
pre-pilot phase.

After gaining knowledge to make the pilot study operational, $\mathbf
{x}$ and $\mathbf{y_x}$ are obtained from preferential sampling, and
again we note the presence of randomness in choosing $\mathbf{x}$.
Given $\mathbf{x}$ and $\mathbf{y_x}$, the next set of locations,
$\mathbf{d}\equiv(d_1,\ldots,d_m)^{\prime}$, need to be chosen, for
which there will be a corresponding (based on the latent vector $\mathbf
{s_d}$) $\mathbf{y_d}$. This is the problem considered by~\citetalias
{FG2015}, and their solution follows closely the proposal of~\cite
{Muller1999}. But there is an important difference:~\citeauthor
{Muller1999} considers $\mathbf{d}$ to be a ``design parameter'' that
he clearly treats as non-stochastic (fixed). We would like to
ask~\citetalias{FG2015} the following question: If $\mathbf{x}$ is
considered to be random in the pilot phase, why would $\mathbf{d}$ be
treated as fixed in the main phase of the study?

The authors follow~\citeauthor{Muller1999}'s~\citeyearpar{Muller1999}
proposal closely but, in our opinion, they lose an opportunity to build
a sequential-sampling-design strategy that updates the posterior
distribution for $\theta$ and $s(\cdot)$ through \textit{random}
sampling from the distribution $[\mathbf{d}|\mathbf{x},\mathbf{y_x}]$.
This can be obtained from $[y(\cdot)|s(\cdot)]$ and $[\mathbf{x},\mathbf
{d}|s(\cdot)]=[\mathbf{d}|\mathbf{x},s(\cdot)][\mathbf{x}|s(\cdot)]$.
The second factor in the product is a log-Gaussian Cox process on $A$,
and the first factor is presumably a log-Gaussian Cox process too, but
on $A\backslash\mathbf{x}$. Can the authors comment on our suggestion
that designed randomness be used to obtain $\mathbf{d}$?

The following section makes strong links between~\citetalias{FG2015}'s
proposal and the survey-sampling literature. It also reinforces the
general desire for a component of randomisation in the design.

\section{Spatial Sampling Designs}
\label{section:spatial}

The development in~\citetalias{FG2015} (and in the literature that it
refers to) frames the sampling problem as the selection of a set of
points on a grid that ``covers'' the region of interest, $A$. From this
perspective, it is a special case of a finite-population sampling
problem, with the $N$ grid points defining the population, and with two
random quantities defined on these points. The first is the
sample-selection process that results in $\mathbf{x}$ (from the pilot
study) and $\mathbf{d}$ (from the main survey). The second is the
latent process $s(\cdot)$, from which ``noisy'' observations $\mathbf
{y_x}$ and $\mathbf{y_d}$ are taken at $\mathbf{x}$ and $\mathbf{d}$,
respectively.

The aim of a spatial sampling design should be to specify a suitable
procedure for making a draw from the distribution of $\mathbf{d}$ given
$\mathbf{x}$ and $\mathbf{y_x}$; see our discussion in Section~\ref
{sect:FixedRandom}. What is meant by ``suitable'' depends crucially on
the target of inference for the sampling exercise;~\citetalias{FG2015}
make their target $s(\cdot)$ and to a lesser extent $\theta$, and they
assume that ``suitability'' can be characterised through a utility
function $u$. Their optimal sample $\mathbf{d}$ is then the set of
(presumably so far unsampled) grid points that maximise the expected
value of this utility, where recall that the expectation is with
respect to the joint distribution of $\theta$ and $\mathbf{s_d}$
conditional on $\mathbf{x}$ and $\mathbf{y_x}$.

The authors' optimal-design procedure is explicitly model-based.
Furthermore, the fact that selection of $\mathbf{d}$ depends on a
log-Gaussian Cox process with intensity function that is a function of
$s(\cdot)$ means that the sampling design is informative \citep[Section
1.4]{ChambersClark2012}, which~\cite{Diggle2010} and~\citetalias
{FG2015} refer to as preferential sampling. That is, one cannot treat
the realised value of $\mathbf{d}$ as ancillary when using the combined
pilot-study and main-survey data to make inferences. There is a well
developed theory in the sample-survey literature for the analysis of a
sample collected via informative sampling; see~\cite
{ChambersSteel2012}. From this perspective, the use of a log-Gaussian
Cox process as a model for $\mathbf{x}$ is equivalent to Poisson
sampling with inclusion probabilities that depend on the values of
$s(\cdot)$ over the grid defining the finite population. We would like
to draw attention to an extensive literature on this type of sampling
and its implications, including many Bayesian approaches; see~\cite
{Nandram2013}.

More complex informative-sampling methods have also been investigated,
principally in the context of sampling spatially clustered populations;
see~\cite{Rapley2008} for a Bayesian specification and, in the context
of sampling on networks, see~\cite{Thompson1996}. Awareness of this
closely related literature would seem advantageous for further
development of the ideas set out in~\citetalias{FG2015}'s paper.

The main inferential paradigm in survey sampling is design-based
inference, which assumes that $y(\cdot)$ is fixed, and all inference is
relative to the distribution of $\mathbf{d}$. Moreover, the outcome of
the pilot study ($\mathbf{x}$ and $\mathbf{y_x}$) is treated as fixed.
Generally, the survey-sampling approach is based on frequentist
inference about population summary statistics. The inference uses
weights obtained from the randomisation in the design, along with the
population values $y(\cdot)$ over the grid. In the simplest case, these
weights are defined by the inverses of the inclusion probabilities for
each of the elements of $\mathbf{d}$, but more general ``calibrated''
weights are typically preferred; see~\cite{Deville1992}. When
informative (i.e. preferential) sampling is used, design-based
inference, although theoretically still applicable, becomes problematic
in practice. In order to carry out the survey sampling, one has to have
access to the distribution of $\mathbf{d}$, which depends on the latent
process $s(\cdot)$. For design-based inference, one might try replacing
$s(\cdot)$ in the intensity function of the log-Gaussian Cox process
with $z(\cdot)$, a spatial covariate whose value is known for every
point on the grid and which is (hopefully) highly correlated with
$s(\cdot)$. This is the model underpinning size-biased sampling;
see~\cite{Patil1978}.

A model-based approach seems therefore necessary under preferential
sampling, such as assuming a spatial-statistical model for $s(\cdot)$.
However, this does not mean that the basic design-based notions of
randomisation, stratification, and clustering cannot be used in a
preferential-sampling approach, since they are all useful tools that
lead to a better representation of a heterogeneous population. In
particular, what happens when the model~\citetalias{FG2015}-(1)
and~\citetalias{FG2015}-(2) does not adequately describe the spatial
variability in $y(\cdot)$ and $s(\cdot)$? The optimality of $\mathbf
{d}$, and the validity of any consequent inference depends critically
on the appropriateness of this model. This is clearly a weakness,
should the design be for a highly scrutinised environmental study where
scientists are worried not only about the environment but also about
the team of lawyers waiting to litigate! Other problems arise when
there are relatively few such choices of $\mathbf{d}$, irrespective of
the values of $\mathbf{x}$ and $\mathbf{y_x}$, all of whose utilities
are comparable.

We would like to reiterate that some form of randomisation in a design
is always a good idea, because it offers protection against a biased
(unintentional or intentional) choice of sample sites (e.g.~\citealp
{AldworthCressie1999}). Perhaps more importantly, randomisation ensures
that an updated fit of an assumed model for $s(\cdot)$ can be validly
assessed from sample data and that replication-based ideas can be used
for this purpose. And finally, when the parametric model is in doubt,
the presence of randomisation allows the possibility that design-based
inference could be used.

As far as we are aware, there has been no work on ``robustifying''
inference based on data collected via preferential sampling, in order
to make it less sensitive to model misspecification. Perhaps~\citetalias
{FG2015}'s paper will stimulate such investigations. Recent research
reported in~\cite{Welsh2013} may provide an indication of how a robust
preferential-sampling approach might work, with these authors
developing an approach to sampling design that minimises the maximum
prediction error in a \textit{neighbourhood} of an assumed model for
$s(\cdot)$. A related line of research concerns what could be termed as
a composite approach to preferential sampling, where a proportion of
the sampling effort is randomly spread over the spatial grid, with the
rest allocated to a more targeted preferential sampling design. An
important research question here concerns how this allocation might be
determined, based on the information in $\mathbf{x}$ and $\mathbf
{y_x}$; this is discussed further in Section~\ref{section:utility}.

We conclude this section by stating some basic elements of a good
sampling design, be it spatial or not. A good design will stratify to
ensure sampling over a range of levels of factors or a range of values
of covariates. A good design will specify, in advance, inference
thresholds and determine the number of observations per stratum needed
to achieve those thresholds. Such designs create a rational basis for
the inevitable compromise between the cost of the study and the ability
to make scientific inferences from incomplete and noisy data (e.g.
\citealp{Cressie1998};~\citealp{Zidek2000}). Finally, a good design
will involve a component of designed randomness, from which
non-parametric, design-based inference is also possible, should the
model-based assumptions be in doubt.

\section{Utility Functions}
\label{section:utility}
The process $s(\cdot)$ and the behaviour of the observations $\mathbf
{y_x}$ depend on parameters, which are denoted as $\theta$. If $\theta$
were known, then $\left[\mathbf{x},\mathbf{y_x},s(\cdot)|\theta\right
]=\left[\mathbf{y_x}|\mathbf{x},s(\cdot),\theta\right]\left[\mathbf
{x},s(\cdot)|\theta\right]$, and the predictive distribution is
\begin{equation}
\label{eq:pred-dist}
\left[s(\cdot)|\mathbf{x},\mathbf{y_x},\theta\right]=\left[\mathbf
{y_x}|\mathbf{x},s(\cdot),\theta\right]\left[\mathbf{x},s(\cdot)|\theta
\right]/\left[\mathbf{x},\mathbf{y_x}|\theta\right].
\end{equation}
Using the terminology of~\cite{CressieWikle2011}, an \textit{empirical
hierarchical model (EHM)} results if an estimate $\widehat{\theta}$ is
used in place of $\theta$ in (\ref{eq:pred-dist}), and inference on
$s(\cdot)$ is then based on the \textit{empirical} predictive distribution,
\begin{equation}
\label{eq:EPD}
[s(\cdot)|\mathbf{x},\mathbf{y_x},\widehat{\theta}]=[\mathbf
{y_x}|\mathbf{x},s(\cdot),\widehat{\theta}][\mathbf{x},s(\cdot)|\widehat
{\theta}]/[\mathbf{x},\mathbf{y_x}|\widehat{\theta}].
\end{equation}
This EHM set-up is what~\cite{Diggle2010} use, and they address the
importance of making $\widehat{\theta}$ a function of both $\mathbf{x}$
and $\mathbf{y_x}$.

If there is uncertainty in $\theta$ that can be expressed in terms of a
prior probability distribution $\left[\theta\right]$, then a \textit
{Bayesian hierarchical model (BHM)} results. Bayes' Theorem yields the
posterior distribution,
\begin{equation}
\label{eq:PD}
\left[s(\cdot),\theta|\mathbf{x},\mathbf{y_x}\right]=\left[\mathbf
{y_x}|\mathbf{x},s(\cdot),\theta\right]\left[\mathbf{x},s(\cdot)|\theta
\right]\left[\theta\right]/\left[\mathbf{x},\mathbf{y_x}\right].
\end{equation}
For a BHM, the \textit{Bayesian} predictive distribution is the
integral of (\ref{eq:PD}) with respect to $\theta$, namely $\int\left
[s(\cdot),\theta|\mathbf{x},\mathbf{y_x}\right]\mathrm{d}\theta$.

The BHM is coherent in the sense that all inferences emanate from a
well defined joint probability distribution. On the other hand, it
requires specification of a prior $\left[\theta\right]$, and it often
consumes a large amount of computing resources. The EHM represents a
compromise that may achieve computational efficiency.

\cite{Diggle2010} do \textit{not} address optimal spatial design in the
way that~\citetalias{FG2015} do. If one sets about doing it, analogous
to~\citetalias{FG2015}'s approach but within~\citeauthor{Diggle2010}'s
EHM framework, one would modify \eqref{eq:exp-utility} so that the
right-hand side would be the expectation taken over $[\mathbf
{y_d}|\mathbf{x},\mathbf{y_x},\theta]$, and hence one would write the
expected utility as $U(\mathbf{d},\theta;\mathbf{x},\mathbf{y_x})$.
Then the \textit{empirical} utility is $U(\mathbf{d},\widehat{\theta
},\mathbf{x},\mathbf{y_x})$ and, analogous to~\citetalias{FG2015}'s
approach, one would find the EHM-optimal $\mathbf{d}$ by maximising
$U(\mathbf{d},\widehat{\theta};\mathbf{x},\mathbf{y_x})$ with respect
to $\mathbf{d}$. Is there a more principled way to account for $\theta$
(which is considered fixed but unknown) in the EHM framework?

Let $W\equiv\left\{\mathbf{d},\theta,\mathbf{y_d}\right\}$ denote all
the unknowns in~\citetalias{FG2015}'s model. Let $\widehat{W}(\mathbf
{x},\mathbf{y_x})$ be one of many possible decisions about $W$ based on
$\mathbf{x}$ and $\mathbf{y_x}$. Some decisions are better than others,
which can be quantified through a very general utility function that is
bounded above, and which we denote as $\mathcal{U}(W,\widehat{W}(\mathbf
{x},\mathbf{y_x}))$; the utility function, $u$, used by~\citetalias
{FG2015} represents a particular form of the more general $\mathcal{U}$
considered here. Note that $\mathcal{U}$ should account for ``how
much?'' and ``why?'' and could be negative. Obviously, large utilities
are preferred, and it is a consequence of decision theory (e.g.~\citealt
{Berger1985}) that the optimal decision is:
\begin{equation}
\label{eq:opt-dec}
W^{*}(\mathbf{x},\mathbf{y_x})=\mathrm{arg}\,\underset{\widehat
{W}(\mathbf{x},\mathbf{y_x})}{\mathrm{sup}}\left\{\mathrm{E}(\mathcal
{U}(W,\widehat{W}(\mathbf{x},\mathbf{y_x}))|\mathbf{x},\mathbf
{y_x})\right\}.
\end{equation}

Now suppose that the goal is inference on $g(W)$, where $g(\cdot)$ is a
known, scientifically interpretable, possibly multivariate function of
$W$. The answer to this inference problem is found in the predictive
distribution, $\left[g(W)|\mathbf{x},\mathbf{y_x}\right]$. Let $\widehat
{g}$ denote a generic predictor of $g(W)$. The \textit{mean} of the
predictive distribution of $g(W)$, namely $\mathrm{E}(g(W)|\mathbf
{x},\mathbf{y_x})$, is a commonly used predictor, but this is just one
of many possibly summaries of $\left[g(W)|\mathbf{x},\mathbf{y_x}\right]$.

Why use the mean? Because it is straightforward to show that $\mathrm
{E}(g(W)|\mathbf{x},\mathbf{y_x})$ solves (\ref{eq:opt-dec}) when the
utility function is ``negative squared-error,'' $-(\widehat
{g}-g(W))^{\prime}(\widehat{g}-g(W))$. However, a negative
squared-error utility assumes equal consequences for under-estimation
as for over-estimation, which is not appropriate when $g(W)$ represents
extreme events, such as crop failure due to drought.

Notice that we have written the utility as a function of all the
unknowns, $W$, and a decision about all the unknowns, $\widehat{W}$.
This gives us the opportunity to design for making inference on $\mathbf
{d}$ simultaneously with making inference on $\theta$, for example.
Recall from Section~\ref{section:spatial} our discussion of the
composite approach to optimal design. One of the components of $\theta$
might be the derivative of the variogram of $s(\cdot)$ at the origin (a
critical parameter for kriging), which we simultaneously want to infer
along with predicting the hidden spatial process $s(\cdot)$.~\cite
{LaslettMcBratney1990} give a composite spatial design that distributes
sampling locations regularly over $A$ (for inference on $s(\cdot)$)
and, around some of those locations, further locations are chosen very
close together (for inference on $\theta$). Do~\citetalias{FG2015} have
any suggestions as to how an \textit{optimal} composite spatial design
might be obtained under their utility-function approach?

In conclusion, we thank the authors for their stimulating paper, and we
can see a number of very interesting research problems waiting to be solved.

%% Appendix %%
% \appendix
% \section{}\label{}

%% Supplement Material %%
% \begin{supplement}
% \sname{}\label{}
% \stitle{}
% \slink[url]{}
% \stitlepost{\\} %jei doi reikia nukelti i kita eilute
% \sdescription{}
% \end{supplement}

%% References %%
%

% Acknowledgements
% \begin{acknowledgement}
% \end{acknowledgement}

\end{document}